\documentstyle[12pt]{article}
\newcommand{\lra}{\longrightarrow}
\newcommand{\complex}{{\bf C}}

\newcommand{\integer}{{\bf Z}}
\newcommand{\rational}{{\bf Q}}

\newcommand{\cD}{{\cal D}}
\newcommand{\cE}{{\cal E}}

\newcommand{\cO}{{\cal O}}

\newcommand{\Pic}{{\rm Pic}}

\newtheorem{theorem}{Theorem}[subsection]
\newtheorem{example}[theorem]{Example}
\newtheorem{proposition}[theorem]{Proposition}
\newtheorem{prop}[theorem]{Proposition}

\newtheorem{corollary}[theorem]{Corollary}
\newtheorem{remark}[theorem]{Remark}
\newtheorem{definition}[theorem]{Definition}

\newtheorem{conjecture}[theorem]{Conjecture}
\newcommand{\prf}{{\it Proof:} }
\newcommand{\qed}{\hspace*{\fill}\hbox{$\Box$}}

\begin{document}
\title{Conifold Transitions 
and Mirror  
Symmetry for  Calabi-Yau Complete Intersections in Grassmannians.}
\author{{\sc  Victor V. Batyrev}\thanks {
Mathematisches Institut, Eberhard-Karls-Universi\"at  T\"ubingen,  
D-72076 T\"ubingen, Germany, 
{\em email address:} batyrev@bastau.mathematik.uni-tuebingen.de},\and
{\sc Ionu\c t Ciocan-Fontanine }\thanks {
 Department of Mathematics, Oklahoma State University, Stillwater, OK 74078, USA,
{\em email address:} ciocan@math.okstate.edu},
\and{\sc Bumsig Kim}\thanks{
Institut Mittag-Leffler, Aurav\" agen 17, S-182 62, Djursholm, Sweden, 
{\em email addresses:} kimbum@ml.kva.se , bumsig@math.kaist.ac.kr},  
\and 
{\sc and Duco van Straten}\thanks {FB 17, Mathematik, 
Johannes Gutenberg-Universit\"at Mainz, D-55099 Mainz, Germany, 
{\em email address:} straten@mathematik.uni-mainz.de}}

\date{}
\maketitle

\thispagestyle{empty}

\begin{abstract} 
In this paper we show that conifold transitions 
between Calabi-Yau $3$-folds can be used for the
 construction of mirror manifolds 
and for the computation of  the instanton numbers  of rational 
curves on  complete intersection Calabi-Yau $3$-folds in Grassmannians. 
Using a natural degeneration of Grassmannians $G(k,n)$ 
to some  Gorenstein toric Fano varieties $P(k,n)$ with 
conifolds singularities which was recently
 described by Sturmfels, we suggest  an explicit 
mirror construction for Calabi-Yau complete intersections 
$X \subset G(k,n)$ of arbitrary dimension. 
Our mirror construction is consistent with the formula for 
the Lax operator conjectured by Eguchi, Hori and 
Xiong for gravitational quantum cohomology of Grassmannians.

\end{abstract}

\newpage

\tableofcontents

\newpage

\section{Introduction}

One of the simplest ways to connect moduli spaces of 
two Calabi-Yau $3$-folds $X$ and $Y$ is a so called {\em conifold 
transition} that attracted interest of physicists 
several years ago in connection 
with {\em black hole condensation} \cite{S,GMS,CGGK}. 
The idea of the conifold transition goes back to Miles Reid \cite{R}, 
who proposed to connect the  moduli spaces of two Calabi-Yau $3$-folds 
$X$ and  $Y$ by choosing a 
point $x_0$ on the moduli space of complex structures on $X$ corresponding 
to a Calabi-Yau $3$-fold $X_0$ 
whose singularities consist of finitely many nodes. If $Y$ is a small 
resolution of singularities on $X_0$ 
which   replaces  the nodes by a union of ${\bf P}^1$'s with normal 
bundle ${\cal O}(-1) \oplus {\cal O}(-1)$, one obtains another 
smooth Calabi-Yau $3$-fold $Y$. 
Let $p$ be the number of nodes on $X_0$, and let $\alpha$ be the number of
relations between the homology classes of the $p$ vanishing $3$-cycles 
on $X$ shrinking to nodes in  $X_0$. Then the Hodge numbers of $X$ and $Y$ 
are related by the following equations  \cite{C}: 
\[ h^{1,1}(Y) = h^{1,1}(X) + \alpha, \]
\[ h^{2,1}(Y) = h^{2,1}(X) -p + \alpha. \]
The  Hodge numbers of mirrors $X^*$ and $Y^*$ of $X$ and $Y$ 
must satisfy the equations  
\[ h^{1,1}(X) = h^{2,1}(X^*), \; h^{1,1}(X^*) = h^{2,1}(X) \] 
and  
\[ h^{1,1}(Y) = h^{2,1}(Y^*), \; h^{1,1}(Y^*) = h^{2,1}(Y). \]
It is natural to expect that the moduli spaces of 
mirrors $X^*$ and $Y^*$ are again connected in the same 
simplest way: i.e., that  
$X^*$ can be obtained by a small 
resolution of some  Calabi-Yau $3$-fold $Y^*_0$ with $p^*$ nodes and
$\alpha^*$ relations, corresponding 
to a  point $y_0^*$ on the moduli space of complex structures on $Y^*$.
Hence, as suggested in \cite{CGGK,GMS,LS} and \cite{DM1}, 
the conifold transition can be used to find mirrors of $X$, 
provided one knows mirrors $Y^*$  of $Y$. 
For this to work, one then needs
$$p^*=\alpha+\alpha^*=p,  $$
i.e., $X_0$ and $Y_0^*$ have the same 
number of nodes and complementary number of relations between them.  
We remark that even for the simplest  
family of Calabi-Yau $3$-folds, quintic hypersurfaces  in ${\bf P}^4$, 
it is an open problem to determine all possible values of $p$ 
\cite{Straten}.

One of the problems solved in this paper is an explicit geometric 
construction of  mirrors $X^*$ for 
 Calabi-Yau complete intersections 
$3$-folds $X$ in Grassmannians $G(k,n)$
(this was only known for quartics in 
$G(2,4)$, as a particular example of complete intersections 
in projective space \cite{LT}). Our method is based on connecting $X$ via a 
conifold transition to complete intersections $Y$ in a toric manifold.
This manifold is a small crepant desingularization 
$\widehat{P(k,n)}$ of a Gorenstein toric 
Fano variety $P(k,n)$, which in turn is a flat degeneration of $G(k,n)$ in its
Pl\"ucker embedding, constructed by Sturmfels (see \cite{St}, Ch. 11). 
Since one knows how to construct mirrors for Calabi-Yau complete 
intersections in $\widehat{P}(k,n)$ \cite{BS,LB}, 
it remains to find an appropriate specialization of the toric mirrors $Y^*$ 
for $Y$ to conifolds $Y_0^*$ whose small resolutions provide mirrors $X^*$ of 
$X$. The choice of the $1$-parameter 
subfamily of $Y_0^*$ among toric mirrors $Y^*$ 
is determined by the monomial-divisor mirror correspondence and the 
embedding 
$${\bf Z} \cong Pic(P(k,n)) \hookrightarrow Pic(\widehat{P}(k,n)) 
\cong {\bf Z}^{ 1 + (k-1)(n-k-1)}.$$ 

We expect that this method of mirror constructions can be applied 
to all Calabi-Yau $3$-folds whose moduli spaces are connected by
conifolds transitions to the web of Calabi-Yau complete intersections in 
Gorenstein toric Fano varieties. This web has been studied in 
\cite{ACJM,AKMS,BLS} as a generalization of the earlier results 
on  Calabi-Yau complete intersection in products of projective 
spaces and in weighted projective spaces \cite{CDLS,CGH}.

In order to obtain the instanton numbers 
of rational curves on Calabi-Yau complete intersections 
in Grassmannians, we compute  a 
generalized  hypergeometric series $\Phi_X(z)$, 
describing  the monodromy invariant period of $X^*$,  
by specializing a $(1 + (k-1)(n-k-1))$-dimensional generalized  
(Gelfand-Kapranov-Zelevinski) $GKZ$-hypergeometric series for 
the main period of toric mirrors $Y^*$ to 
a single monomial parameter $z$. 
Since $h^{1,1}(X) =1$, the  corresponding   Picard-Fuchs differential 
system for periods of $X^*$ reduces to an ordinary  
differential  equation ${\cal D} \Phi =0$ 
of order $4$  for $\Phi_X(z)$. The Picard-Fuchs 
differential operator $P$ can be computed from the 
recurrent relation satisfied by  the coefficients 
of the series $\Phi_X(z)$.  Applying the same computational 
algorithm as in \cite{BS}, one  computes the instanton numbers of 
rational curves on all possible 
Calabi-Yau complete intersection $3$-folds 
$X \subset G(k,n)$. The numbers of lines and conics on 
these Calabi-Yau $3$-folds have 
been verified  by S.-A. Str{\o}mme 
 using   classical methods and the Schubert package for MapleV. 

Another new ingredient of the present 
paper is the so called {\em Trick with the Factorials}. 
This is a naive form of a 
{\em Lefschetz hyperplane section theorem} in quantum cohomology,
 which goes back to Givental's idea \cite{G2} 
about the relation between solutions of quantum 
$\cD$-module for Fano manifolds $V$  and complete intersections $X \subset V$. 
The validity of this procedure
has been established recently for all homogeneous spaces by B. Kim 
in \cite{K3}. If the Trick with the Factorials 
works for a Fano manifold $V$, one is able to  compute  
the instanton numbers of rational curves on Calabi-Yau complete 
intersections $X \subset V$ without knowing a mirror $X^*$ for $X$, provided 
one knows a special regular solution 
$A_V$ to the quantum $\cD$-module for $V$. 
In the case of Grassmannians we 
conjecture in \ref{flagmirror} that this special solution $A_{G(k,n)}(q)$
to the quantum $\cD$-module determined by the small quantum cohomology 
of $G(k,n)$) 
can be obtained from 
a natural specialization of the $GKZ$-hypergeometric series 
associated with the Gorenstein toric degeneration  $P(k,n)$ of  $G(k,n)$. 
Conjecture \ref{flagmirror} has been  
checked by direct computation  for all Grassmanians containing 
Calabi-Yau $3$-folds 
$X$ as complete intersections. In fact,
there is no essential difficulty in checking  the conjecture 
in each particular case at hand, because such a check involves only 
calculations in the small quantum cohomology ring of 
$G(k,n)$, whose structure is well-known \cite{Ber}. 
This   last result   implies that the   instanton numbers for 
rational curves on $3$-dimensional 
Calabi-Yau complete intersections in Grassmannians are correct 
in all computed cases. 

We remark that our conjecture \ref{flagmirror} on 
the coincidence of $A_{G(k,n)}(q)$ with the 
specialization of the multidimensional generalized $GKZ$-hypergeometric 
series corresponding to the Gorenstein toric Fano variety 
$P(k,n)$ strongly supports the idea  that
Gromov-Witten invariants of $G(k,n)$ and 
complete intersections $X \subset G(k,n)$ 
behave well under flat deformation and conifold transitions.
 
Using  the degeneration of $G(k,n)$ to  $P(k,n)$, we  
propose in arbitrary dimension an explicit construction for mirrors of 
Calabi-Yau complete intersections $X \subset G(k,n)$ 
whose monodromy invariant period 
coincide with  the power series $\Phi_X(z)$  
obtained by applying the Trick with the Factorials 
to $A_{G(k,n)}(q)$.  We  observe that our mirror construction 
is consistent with the formula 
for the Lax operator of Grassmannians conjectured by 
Eguchi, Hori and Xiong in \cite{EHX}.

Many results formulated  in  this paper have been generalized and proved 
in \cite{BCKS} for toric degenerations of 
partial flag manifolds which have been introduced and investigated 
by  N. Gonciulea and V. Lakshmibai in  
\cite{GL0,GL1,GL2}. These results are most easily interpreted in terms of
certain {\em diagrams} associated to a partial flag manifold, generalizing
the one used in \cite{G2} for the case of the complete flag manifold. 

\section{Simplest Examples}  

\subsection{Quartics in $G(2,4)$} 

First we illustrate our method by analyzing a simple case, for which the mirror
construction is already known: the case of  
quartics in $G(2,4)$, the Grassmannian of 
$2$-planes in $\complex ^4$ \cite{BS,LT}. The 
Pl\"{u}cker embedding realizes
the Grassmannian $G(2,4)$ as a nonsingular quadric in ${\bf P}^5$, defined 
by the homogeneous equation: 
\[ z_{12} z_{34} - z_{13}z_{24} + z_{14}z_{23} = 0,   \]
where $z_{ij}$ $(1\leq i<j \leq 4)$ are homogeneous 
coordinates on ${\bf P}^5$. Let  $P(2,4) \subset {\bf P}^5$ 
be the $4$-dimensional Gorenstein toric Fano variety defined by 
the quadratic equation 
\[ z_{13}z_{24} = z_{14}z_{23}. \] 
Denote by  $X$ the intersection of $G(2,4)$ with a generic 
hypersurface $H$ of degree $4$ in ${\bf P}^5$, so that $X$ is   
a nonsingular Calabi-Yau hypersurface in $G(2,4)$. 
Its topological invariants are
$h^{1,1}(X) = 1$, $h^{2,1}(X) = 89$, and $\chi(X)= -176$. Denote by 
$X_0$ the  intersection of $P(2,4)$ with a generic hypersurface $H$ of 
degree $4$ in ${\bf P}^5$. Then $X_0$ is a Calabi-Yau $3$-fold with 
$4$ nodes which are the intersection points of $H$ with 
the line $l \subset P(2,4)$ 
of conifold singularities. Considering $X_0$ as a deformation 
of $X$, it follows from general formulas proved in \ref{formul-2} 
that the homology classes of the vanishing $3$-cycles on $X$ shrinking 
to $4$ nodes in $X_0$ satisfy a single relation. Denote by 
$Y$ a simultaneous small resolution of all $4$ nodes. One obtains this 
resolution by restriction of a small toric resolution of 
singularities in $P(2,4)$: $\rho\,: \; \widehat{P}(2,4) \rightarrow P(2,4)$. 
The smooth toric variety $\widehat{P}(2,4)$ is a toric ${\bf P}^3$-bundle 
over ${\bf P}^1$: 
\[ \widehat{P}(2,4) = {\bf P}_{{\bf P}^1}({\cal O} \oplus {\cal O} 
\oplus {\cal O}(1) \oplus {\cal O}(1)) \]
and the morphism $\rho$ contracts a $1$-parameter family 
of sections  of this ${\bf P}^3$-bundle with 
the normal bundle ${\cal O} \oplus {\cal O}(-1) \oplus {\cal O}(-1)$. 
A smooth Calabi-Yau hypersurface $Y \subset \widehat{P}(2,4)$ has a natural 
$K3$-fibration over ${\bf P}^1$ and the following 
topological invariants: $\chi(Y) = -168$, $h^{1,1}(Y) = 2$, and 
$h^{2,1}(Y)=86$.

The Gorenstein toric Fano variety $P(2,4)$ can be 
described by a $4$-dimensional fan $\Sigma(2,4) \subset {\bf R}^4$ consisting 
of cones over the faces of a $4$-dimensional reflexive polyhedron 
$\Delta(2,4)$ with $6$ vertices:
\[ u_{1,0} := f_{1,1}, \; u_{2,0} = f_{2,1} - f_{1,1}, \; 
u_{2,1}: = f_{2,2} - f_{1,2},  \] 
\[  v_{2,2} := - f_{2,2}, \;
 v_{2,1}: = f_{2,2} - f_{2,1}, \; 
v_{1,1} = f_{1,2} - f_{1,1},   \]  
where $\{ f_{1,1}, f_{1,2}, f_{2,1}, f_{2,2} \}$ is a basis 
of the lattice ${\bf Z}^4 \subset {\bf R}^4$.

The regular fan $\widehat{\Sigma}(2,4)$ defining the smooth projective 
toric variety $\widehat{P}(2,4)$ 
is obtained by a subdivision  of $\Sigma(2,4)$. 
The combinatorial structure of $\widehat{\Sigma}(2,4)$ is defined by the 
following primitive collections (see notations in \cite{Ba1}): 
\[ {\cal R} = \{ u_{1,0}, v_{1,1}, u_{2,1}, v_{2,2} \} , \; 
 {\cal C}_{1,1} = \{ v_{2,1}, u_{2,0} \}. \]
The fan $\widehat{\Sigma}(2,4)$ contains $8$ cones of dimension $4$, obtained
by deleting one vector from each primitive collection. 
The primitive relations corresponding to ${\cal R}_0$ and ${\cal C}_{1,1}$ 
are 
\[  u_{1,0} +  v_{1,1} +  u_{2,1} + v_{2,2} = 0 \]
and 
\[  v_{2,1} +  u_{2,0} = v_{1,1} + u_{2,1}. \]

Let ${\bf P}_{\Delta(2,4)}$ be the Gorenstein toric Fano variety 
associated with the reflexive polyhedron $\Delta(2,4)$.
By the toric method of \cite{Ba2}, the mirror $Y^*$ of $Y$ can be obtained as 
a crepant desingularization of the closure in ${\bf P}_{\Delta(2,4)}$
of an affine hypersurface    
$Z_f$ with the equation 
\[ f(X) = -1 + a_1X_1 + a_2X_2 + a_3 X_3 + a_4 X_4 + a_5(X_1X_2X_3)^{-1} + 
a_6(X_4^{-1}X_1X_2),  \]
where $a_1, \ldots, a_6$ are some complex numbers 
(the Newton polyhedron of $f$ is isomorphic 
to  $\Delta(2,4)$). 
We choose a subfamily of Laurent polynomials $f_0$ with  coefficients 
$\{a_i \}$ satisfying an  additional monomial equation 
\[ a_1 a_2 = a_4 a_6. \]  
 The affine Calabi-Yau hypersurfaces $Z_{f_0}$ of this subfamily are 
not $\Delta(2,4)$-regular anymore, because the  closures 
$\overline{Z}_{f_0}$ in ${\bf P}_{\Delta(2,4)}$ have  a singular 
intersection with the stratum 
$T_{\Theta} \subset {\bf P}_{\Delta(2,4)}$ corresponding 
to the face 
\[ \Theta = {\rm Conv}\left\{  (1,0,0,0),\, (0,1,0,0), \, 
(0,0,0,1),\,(1,1,0,-1) \right\}. \]
Without loss of generality, we can assume that $a_1=a_2=a_3=a_4 =1$
(this condition can be satisfied using the action of $({\bf C}^*)^4$ on 
$X_1, \ldots, X_4$). Thus  we obtain a $2$-parameter family of Laurent 
polynomials defining $Z_f$:  
\[ f(X) = -1 + X_1 + X_2 + X_3 + X_4 + a_5(X_1X_2X_3)^{-1} + 
a_6(X_4^{-1}X_1X_2), \]
and a $1$-parameter subfamily of Laurent polynomials   
\[ f_0(X) = -1 + X_1 + X_2 + X_3 + X_4 + a_5(X_1X_2X_3)^{-1} + 
(X_4^{-1}X_1X_2) \]
defining $Z_{f_0}$. 
The monodromy invariant period $\Phi$ 
of the toric hypersurface $Z_f$ can be computed by the  residue theorem: 
\[ \Phi_X(a_5,a_6) = \frac{1}{(2\pi i )^4} \int_{\gamma} \frac{1}{(-f)} 
\frac{dX_1}{X_1} \wedge \cdots \wedge \frac{dX_4}{X_4}. \]
By this method, we obtain the generalized 
hypergeometric series corresponding to $f(X)$: 
\[ \Phi_X(a_5,a_6) = \sum_{k, l \geq 0} \frac{(4k + 4l)!}{(k!)^2 (l!)^2 
((k+l)!)^2} a_5^{k+l} a_6^l. \]
By the substitution  $a_6 =1$ $(a_1a_2 = a_4a_6)$ 
and $a_5 = z$, we obtain the series corresponding to 
the $1$-parameter family of Laurent polynomials 
$f_0$: 
\[ \Phi_X(z) = \sum_{m \geq 0} \frac{(4m)!}{(m!)^2} 
\left(\sum_{k+l = m} \frac{1}{(k!)^2 (l!)^2} \right)z^m. \]
Using the identity
\[ \sum_{k+l =m} \frac{(m!)^2}{(k!)^2 (l!)^2} = { 2m \choose m}, \]
we transform $\Phi_X(z)$ to the form 
\[ \Phi_X(z) =     \sum_{m \geq 0} \frac{(4m)!(2m)!}{(m!)^6} z^m. \]
This is a well-known series, satisfying a Picard-Fuchs 
differential equation 
\[ \left(  D^4 - 16z(2 D +1)^2(4 D + 1)(4 D + 3)\right) 
\Phi_X(z) = 0,\;\;  
 D = z\frac{\partial}{\partial z}, \] 
predicting the instanton numbers of rational curves on $X$ (cf. \cite{LT}). 
The correctness of these numbers
 now follows from the work of Givental, \cite{G1}.

\subsection{Complete intersections of type  $(1,1,3)$ in $G(2,5)$}
\vskip 10pt
Let  $z_{ij}$ $(1 \leq i < j \leq 5$ ) be  homogeneous 
coordinates on  the projective space ${\bf P}^9$.   
The Grassmannian $G(2,5)$ of $2$-planes in $\complex^5$ 
can be identified with the subvariety 
in ${\bf P}^9$ 
defined by the quadratic equations: 
\[ z_{23}z_{45} -  z_{24}z_{35} +  z_{25}z_{34}= 0,  \]     
\[ z_{13}z_{45} -  z_{14}z_{35} +  z_{15}z_{34}= 0,  \]     
\[ z_{12}z_{45} -  z_{14}z_{35} +  z_{15}z_{34}= 0,  \]     
\[ z_{12}z_{35} -  z_{13}z_{25} +  z_{15}z_{23}= 0,  \]     
\[ z_{12}z_{34} -  z_{13}z_{24} +  z_{14}z_{23}= 0.  \]     
We associate with $G(2,5)$ a  $6$-dimensional Gorenstein 
toric Fano variety $P(2,5) \subset {\bf P}^9$ defined by the 
equations 
\[ z_{24}z_{35} =  z_{25}z_{34}, \;   z_{14}z_{35} = z_{15}z_{34},\;       
z_{14}z_{35} = z_{15}z_{34},  \]     
\[ z_{13}z_{25} = z_{15}z_{23},  \; z_{13}z_{24} = z_{14}z_{23}.  \]
The following statement is due to Sturmfels (see \cite{St}, Example 11.9 and 
Proposition 11.10): 

\begin{prop}
The Gorenstein toric Fano variety 
$P(2,5)$ is a degeneration of the Grassmannian $G(2,5)$, i.e., 
$P(2,5)$ is the special fibre of a 
flat family whose generic fibre is $G(2,5)$. \qed
\end{prop}

The toric variety $P(2,5)$ can be described by a fan 
$\Sigma(2,5) \subset {\bf R}^6$ consisting of cones over the faces of 
a $6$-dimensional reflexive polyhedron $\Delta(2,5)$ with $9$ vertices 
\[ u_{1,0} := f_{1,1}\;   u_{2,i}: = f_{2,i+1} - f_{1,i+1}, \,\; i = 0, 1, 2,\]
\[ v_{2,3} := - f_{2,3}, \;  v_{i,j}: = f_{i,j+1} - f_{i,j}, \, \; 
i = 1, 2,\; j =1,2,  \]
where $\{ f_{1,1}, f_{1,2}, f_{1,3}, f_{2,1}, f_{2,2}, f_{2,3} \}$
is a basis of the lattice ${\bf Z}^6 \subset {\bf R}^6$.

There exists a subdivison  
of the fan $\Sigma(2,5)$ into a regular 
fan $\widehat{\Sigma}(2,5)$ defined by the primitive collections:
\[ {\cal R} = \{ u_{1,0}, v_{1,1}, v_{1,2}, u_{2,2}, v_{2,3} \}, \]
\[ {\cal C}_{1,1} = \{ u_{2,0}, v_{2,1} \},  \; 
 {\cal C}_{1,2} = \{ u_{2,1}, v_{2,2} \},  \]
i.e., $\widehat{\Sigma}(2,5)$ 
contains exactly $20$ cones of dimension $6$ generated 
by the $6$-element sets obtained by taking all but one of the 
vectors from each
primitive collection. 
The primitive relations corresponding to 
${\cal R}$, ${\cal C}_{1,1}$ and ${\cal C}_{1,2}$ are 
\[ u_{1,0} +  v_{1,1} +  v_{1,2} +  u_{2,2} +  v_{2,3} = 0, \]
\[  u_{2,0} + v_{2,1} = v_{1,1} + u_{2,1},  \; \; 
  u_{2,1} +  v_{2,2} = v_{1,2} + u_{2,2}. \]  

Denote by $\widehat{P}(2,5)$ the smooth toric variety associated with 
the fan $\widehat{\Sigma}(2,5)$. It is easy to check that $P(2,5)$ is a 
Gorenstein toric Fano variety and $\widehat{P}(2,5)$ is a small 
crepant resolution 
of singularities of $P(2,5)$. The toric manifold 
$\widehat{P}(2,5)$ has nonnegative 
first Chern class and it can be identified with a 
toric bundle over ${\bf P}^1$ with the $5$-dimensional fiber 
\[ F: = {\bf P}_{{\bf P}^1} ({\cal O} \oplus {\cal O} \oplus {\cal O} \oplus 
{\cal O}(1) \oplus {\cal O}(1) ) \]

There is another description of $P(2,5)$. We remark that 
variables $z_{12}$ and $z_{45}$ do not appear in the equations for $P(2,5)$. 
Thus   $P(2,5)$ is  a cone over a Gorenstein $4$-dimensional toric Fano 
variety
\[ P'(2,5) : = P(2,5) \cap \{ z_{12} = z_{45} = 0\} \subset {\bf P}^7. \]
We can describe $P'(2,5)$ by a $4$-dimensional fan $\Sigma'(2,5)$ 
consisting of cones over a $4$-dimensional reflexive polyhedron 
$\Delta'(2,5)$ with $7$ vertices
\[ e_1 = (1,0, 0,0), \; e_2= (0,1,0,0),\; e_3 = (-1,-1, 0,0), \]   
\[ e_4 = (0,0, 1,0), \; e_5= (0,0,0,1),\;  e_6 = (0,0,-1,-1), \]  
\[ e_7 = (1,1,1,1). \]
The only singularities of $P'(2,5)$ are 
nodes along two lines $l_1, l_2 \in P'(2,5) \subset {\bf P}^7$ 
corresponding to the 
$3$-dimensional cones 
$$ \sigma_1 = {\bf R}_{\geq 0} < e_1, e_2, e_6, e_7 >\;\; 
\mbox{\rm and } \;\;\sigma_2 = {\bf R}_{\geq 0} 
< e_4, e_5, e_3, e_7 >$$ 
in $\Sigma'(2,5)$. 
Subdividing each of these cones into the union of $2$ simplicial ones, we 
obtain a small crepant resolution $\widehat{P'}(2,5)$ of singularities 
of $P'(2,5)$. The smooth toric  $4$-fold $\widehat{P'}(2,5)$ can be identified 
with the blow up of a point on ${\bf P}^2 \times {\bf P}^2$. 

Let $X = X_{1,1,3} \subset G(2,5)$ be a  
smooth $3$-dimensional Calabi-Yau complete intersection 
of $3$ hypersurfaces of degrees $1$, $1$ and $3$ in ${\bf P}^9$ 
with $G(2,5)$. One can compute $h^{1,1}(X) = 1$,   $h^{2,1}(X) = 76$,
and $\chi(X) = -150$. Now let $X_0$ be the intersection of $P'(2,5)$ with 
a generic hypersurface $H \subset {\bf P}^7$ of degree $3$. Then 
$X_0$ is a deformation of $X$, 
having $6$ nodes obtained from the intersections
$H \cap l_1$ and $H \cap l_2$. The $3$ nodes on each intersection   
$H \cap l_i$ $(i =1,2)$ are  described by $3$ vanishing $3$-cycles on 
$X$, satisfying a single 
linear relation. Resolving singularities of $X_0$, we obtain 
another smooth Calabi-Yau $3$-fold $Y$ with 
\[ h^{1,1}(Y) = h^{1,1}(X) + 2 = 3, \;\; 
h^{2,1}(Y) = h^{2,1}(X) + 2 - 6 = 72. \]      

The mirror $Y^*$ of the Calabi-Yau $3$-fold  $Y$ 
can be obtained by the toric construction
\cite{Ba2}. The Calabi-Yau $3$-fold $Y^*$ is  a toric 
desingularization  $\widehat{Z}_f$ of a $\Delta'(2,5)$-compactification of a  
generic hypersurface $Z_f$ in $({\bf C}^*)^4$ 
defined by a  Laurent polynomial $f(X)$ with the Newton polyhedron 
$\Delta'(2,5)$: 
\[ f(X) = -1 + a_1 X_1 + a_2 X_2 + a_3 (X_1X_2)^{-1} + a_4 X_3 + a_5X_4+ \] 
\[  + a_6(X_3X_4)^{-1} + a_7X_1X_2X_3X_4. \]
As it was shown in \cite{Ba2}, one has   
$h^{1,1}(\widehat{Z}_f) = h^{2,1}(Y) = 72$ and 
$h^{2,1}(\widehat{Z}_f) = h^{1,1}(Y) = 3$. 

We identify the mirror $X^*$ of $X$ 
with a desingularization $\widehat{Z}_{f_0}$ of 
a $\Delta'(2,5)$-compactification $\overline{Z}_{f_0}$ of a  
generic hypersurface $Z_{f_0}$ in $({\bf C}^*)^4$ 
defined by Laurent polynomials $f_0$ whose 
coefficients $\{ a_i\} $ satisfy two additional 
monomial equations
$$ a_1a_2 = a_6a_7\ \ \mbox{\rm and} \ \ a_4 a_5 = a_3 a_7. $$
Without loss of generality, we can put $a_1 = a_2 = a_4 = a_7$. So one 
obtains
\[ f(X) = -1 + X_1 + X_2 + a_3 (X_1X_2)^{-1} + X_3 + a_5X_4 \] 
\[  + a_6(X_3X_4)^{-1} + X_1X_2X_3X_4 \]
and 
\[ f_0(X) = -1 + X_1 + X_2 + a_3 (X_1X_2)^{-1} + X_3 + a_3X_4 \] 
\[  + (X_3X_4)^{-1} + X_1X_2X_3X_4. \]
It is easy to see that the Laurent polynomial $f_0$ 
is not $\Delta'(2,5)$-regular (this regularity fails 
exactly for two $2$-dimensional faces $\Theta_1 
:= {\rm Conv}(e_1,e_2,e_6,e_7)$ 
and  $\Theta_2:= {\rm Conv}(e_4,e_5,e_3,e_7)$ 
of $\Delta'(2,5)$ (see definition of $\Delta$-regularity in  
\cite{Ba2}). 
The $4$-dimensional 
Gorenstein toric Fano variety ${\bf P}_{\Delta'(2,5)}$ associated with 
the reflexive polyhedron $\Delta'(2,5)$-closure  has singularities 
of type $A_2$ along of the $2$-dimensional strata $T_{\Theta_1}$ and 
$T_{\Theta_2}$. The projective hypersurfaces  $\overline{Z}_{f_0} \subset 
{\bf P}_{\Delta'(2,5)}$ 
defined by the equation $f_0 =0$ have non-transversal intersections with 
$T_{\Theta_1}$ and $T_{\Theta_2}$ (each  intersection is a  union of 
two rational curves with a single normal crossing point).  After toric 
resolution of $A_2$-singularities along $T_{\Theta_i}$ 
on ${\bf P}_{\Delta'(2,5)}$, we obtain $3$ new $2$-dimensional 
strata over each $T_{\Theta_i}$. This shows that 
we cannot resolve all singularities of $\overline{Z}_{f_0}$ by a toric 
resolution of  singularities on the ambient toric variety 
${\bf P}_{\Delta'(2,5)}$. 
Let $Y_0^* := \widehat{Z}_{f_0}$ be the pullback of $\overline{Z}_{f_0}$ under 
a $MPCP$-desingularization 
$$\rho\, : \, \widehat{{\bf P}}_{\Delta'(2,5)} \rightarrow 
{\bf P}_{\Delta'(2,5)}.$$ 
Then $Y_0^*$ is a Calabi-Yau 3-fold with
$3 + 3 = 6$ nodes obtained as singular points of intersections 
of $Y_0^*$ with 
the $6$ strata of dimension $2$  in $\widehat{{\bf P}}_{\Delta'(2,5)}$
over $T_{\Theta_1},\; T_{\Theta_2} \subset  \widehat{{\bf P}}_{\Delta'(2,5)}$. 
One can show that the 
vanishing $3$-cycles associated with the $3$ nodes over each 
$T_{\Theta_i}$ $(i=1,2)$ satisfy $2$ linear relations (see \ref{formul-2}). 
If $X^*$ denotes a small resolution of 
these $6$ nodes on $Y_0^*$, then  
\[ h^{1,1}(X^*) = h^{1,1}(\widehat{Z}_{f}) + 4 = 76 \]
and 
\[ h^{2,1}(X^*) = h^{2,1}(\widehat{Z}_f) + 4 - 6 = 1. \]      
Thus the Hodge numbers of $X^*$ and $X$ satisfy the mirror duality.

Finally, we  explain the computation of the instanton numbers 
of rational curves 
of degree $m$ in  the case of 
Calabi-Yau complete intersections $X$ of type $(1,1,3)$ 
in $G(2,5)$. As shown in \cite{BS}, one obtains 
the following monodromy invariant period for $Z_f$:
$$\Phi(a_3,a_5,a_6) = \sum_{k,l,n \geq 0} 
\frac{(3k+3l+3n)!}{(k!)^2(n!)^2l!(k+l)!(l+n)!}a_3^{k+l}a_5^n a_6^{n+l}.$$ 
By the substitution $a_3 = a_5 =z$ and $a_6 = 1$, we obtain the monodromy 
invariant period for $Z_{f_0}$: 
\[ \Phi_X(z) = \left( \sum_{k + l +n = m} 
\frac{(3m)!}{(k!)^2(n!)^2l!(k+l)!(l+n)!} \right) z^m. $$
It remains to apply to the series $\Phi_X(z)$ the general algorithm 
from \cite{BS} (see 6.2 and 7.1 for details, and the instanton numbers). 
\vskip 10pt

\section{Toric Degenerations of Grassmannians}

In this section we review without proof some results, which we prove 
for arbitrary partial flag manifolds in \cite{BCKS}. 

\subsection{The toric variety $P(k,n)$ and its singular locus}

 Let $G(k,n)$ be the Grassmannian of 
$k$-dimensional ${\bf C}$-vector subspaces
in a $n$-dimensional complex vector space ($k < n$). 
Denote by  
$$X_{i,j} \;\; i = 1, \ldots, k, \; j = 1, \ldots, n-k$$ 
$k(n-k)$ 
independent variables. We denote by $T(k,n)$  the algebraic torus 
$ {\rm Spec}\, {\bf C}[X_{i,j},X_{i,j}^{-1} ]\cong 
({\bf C}^*)^{k(n-k)}$ of dimension $k(n-k)$. 
We put $N(k,n): = {\bf Z}^{k(n-k)} $ to be a free abelian group of rank 
$k(n-k)$ with a fixed   ${\bf Z}$-basis $f_{i,j}$ 
$(i = 1, \ldots, k, \; j = 1, \ldots, n - k)$. 
Define the set of $2(k-1)(n-k-1) + 
n$ elements 
in $N(k,n)$ as follows:
\[ u_{1,0} := f_{1,1},\;  u_{i,j}: = f_{i,j+1} - f_{i-1,j+1}, 
\, \; i = 2, \ldots, k,\; 
j =0, \ldots, n-k-1\; \]
\[  v_{k,n-k}: = - f_{k,n-k}, \; 
 v_{i,j}: = f_{i,j+1} - f_{i,j}, \,\; i = 1, \ldots, k,\; 
j =1, \ldots, n-k-1.  \]

We set $N(k,n)_{\bf R} = N(k,n) \otimes {\bf R}$. 

\begin{definition}
{\rm Define a convex polyhedron 
$\Delta(k,n) \subset  N(k,n)_{\bf R}$ as the convex hull 
of all  lattice points $\{ u_{i,j}, v_{i',j'} \}$. 
We set $\Sigma(k,n) \subset N(k,n)_{\bf R}$ to be the  fan over all 
proper faces of the polyhedron $\Delta(k,n)$.} 
\label{polyh}
\end{definition} 

\begin{definition}
{\rm Define $P(k,n)$ to be the toric variety associated with  
the fan $\Sigma(k,n)$. } 
\end{definition}

\begin{theorem} 
The polyhedron $\Delta(k,n)$ is reflexive. In particular, 
$P(k,n)$ is a Gorenstein toric Fano variety. 
\end{theorem} 

\begin{definition} 
{\rm Let $\widehat{\Sigma}(k,n)$  be a  complete regular fan 
whose $1$-dimensional cones are generated by the lattice vectors  
$\{  u_{i,j}, v_{l,m}\} $ and whose combinatorics  
is defined by the following $1 + (k-1)(n-k-1)$ primitive  collections: 
\[ {\cal R}_0: = \{ u_{1,0}, v_{1,1}, v_{1,2},  \ldots, v_{1,{n-k-1}}, 
u_{2, n-k-1}, u_{3, n-k -1}, \ldots, u_{k,n-k-1}, v_{k,n-k} \}, \]
\[ {\cal C}_{i,j} = \{ u_{k+1-i, j-1}, v_{k+1-i, j},\;  
\; i=1, \ldots, k-1,\;  j =1, \ldots, n-k-1 \}.\]
In particular, the fan $\widehat{\Sigma}(k,n)$ 
consists of $n2^{(k-1)(n-k-1)}$ 
cones of dimension $k(n-k)$.
}
\end{definition}

\begin{remark}
{\rm We notice that the lattice vectors $u_{i,j}$ and $v_{l,m}$ satisfy the 
following $1+ (k-1)(n-k-1)$ independent primitive relations: 
\[  u_{1,0} +  v_{1,1} +  \cdots +  v_{1,{n-k-1}} +  
u_{2, n-k-1} +  \cdots +  u_{k,n-k-1} +  v_{k,n-k} = 0, \]
\[  u_{k+1-i, j-1} +  v_{k+1-i, j} = u_{k+1-i, j} + v_{k-i,j}, \]  
\[ i=1, \ldots, k-1,\;  j =1, \ldots, n-k-1 .\]
According to  Theorem 4.3 in \cite{Ba1}, the toric variety 
 $\widehat{\Sigma}(k,n)$ can be obtained as $ (k-1)(n-k-1)$-times 
iterated toric bundle over ${\bf P}^1$'s: we start with ${\bf P}^{n-1}$ 
and  construct on  each step a toric bundle over ${\bf P}^1$ 
whose fiber is the toric variety constructed in the previous step.  
At each stage of this process, we obtain a smooth projective toric variety 
with the nonnegative first Chern class which is divisible by $n$.
In particular we obtain that the smooth projective toric 
variety $\widehat{P}(k,n)$ defined by the fan $\widehat{\Sigma}(k,n)$
has Picard number $1 + (k-1)(n-k-1)$. Moreover, the first Chern class 
$\widehat{c}_1(k,n)$ of 
$\widehat{P}(k,n)$ is nonnegative and it is divisible by $n$ 
in ${\rm Pic}(\widehat{P}(k,n))$.
 }
\end{remark}

\begin{definition}
We denote by  
$\widehat{P}(k,n)$ $( 1 \leq i \leq k-1, \; 
1 \leq j \leq n- k-1)$  $(k-1)(n-k-1)$ codimension-$2$ subvarieties 
of $\widehat{P}(k,n)$ corresponding to the $2$-dimensional cones  $
\sigma_{ij} \in  \widehat{\Sigma}(k,n)$:  
\[ \sigma_{ij} = {\bf R}_{\geq 0} <u_{k+1-i, j-1}, 
 v_{k+1-i, j} >. \]  
\end{definition}

\begin{theorem}
The small contraction $\rho\, : \, \widehat{P}(k,n) \rightarrow {P}(k,n)$ 
defined by the semi-ample anticanonical divisor 
on  $\widehat{P}(k,n)$  
contracts smooth toric varieties 
$\widehat{W}_{i,j}$ to codimension-$3$ toric subvarieties $W_{i,j} 
\subset {P}(k,n)$ whose open strata consist of conifold singularities, 
i.e., 
singularities whose $3$-dimensional cross-sections are isolated nondegenerate 
quadratic singularities (nodes, ordinary double points).  
\label{sing-l}
\end{theorem}

The proof of a generalized version of \ref{sing-l} for arbitrary 
partial flag manifolds is contained in \cite{BCKS}(Th. 3.1.4).

\subsection{The flat degeneration of 
$G(k,n)$  to $P(k,n)$}

\begin{definition}
{\rm Denote by $A(k,n)$ the set of all sequences of  integers
\[ a= (a_1, a_2, \ldots, a_k)  \in {\bf Z}^k \] 
satisfying the condition
\[  1 \leq a_1 <  a_2 <  \cdots < a_k \leq n. \] 
We consider $A(k,n)$ as a partially ordered set 
with the  following natural partial  order: 
\[ a = (a_1, \ldots, a_k) \prec a'= (a_1', \ldots, a_k') \]
if and only if $a_i \leq a_i'$ for all $i =1, \ldots, k$.  
We set 
\[ \min{(a,a')} : = (\min{(a_1,a_1')}, \ldots, \min{(a_k,a_k')}) \]
and 
\[  \max{(a,a')} : = (\max{(a_1,a_1')}, \ldots, \max{(a_k,a_k')}). \]  } 
\end{definition}

\begin{theorem} 
There exists a natural one-to-one correspondence between 
faces of codimension $1$ of the polyhedron $\Delta(k,n)$ 
and 
elements of $A(k,n)$.  
\end{theorem} 

\noindent
{\em Proof.} See \cite{BCKS} (Th. 2.2.3). 

\begin{theorem}
The first Chern class 
of the Gorenstein 
toric Fano variety  ${P}(k,n)$ is equal to  $n[H]$, 
where $[H]$ is the class of the ample generator of  
${\rm Pic}({P}(k,n)) \cong {\bf Z}$. Moreover, there 
exists a natural one-to-one correspondence 
between the elements of the monomials basis of 
\[ H^0({P}(k,n), {\cal O}(H)) \]
and elements of 
$A(k,n)$.  
In particular, 
\[ {\rm dim} \, H^0({P}(k,n), {\cal O}(H)) = {  n \choose k }. \]
\label{gl-sections}
\end{theorem}

\noindent
{\em Proof.} See \cite{BCKS} (Prop. 3.2.5).  

\begin{theorem} 
The ample line bundle ${\cal O}(H)$ on $P(k,n)$ 
defines a projective embedding into the
projective space ${\bf P}^{{ n \choose k}-1}$ 
whose homogeneous coordinates $z_a$ are naturally indexed 
by elements $a \in A(k,n)$. Moreover, 
the image of  $P(k,n)$ in  ${\bf P}^{{ n \choose k}-1}$ 
is defined by the quadratic homogeneous binomial equations
\[  z_a z_{a'} - z_{min(a,a')}z_{max(a,a')}  \]
for all pairs $(a,a')$ of non-comparable elements $a, a' \in 
A(k,n)$. 
\end{theorem} 

\noindent
{\em Proof.} See \cite{BCKS} (Th. 3.2.13).

\begin{example}
{\rm The following ${ n \choose 4 }$ quadratic 
equations in homogeneous coordinates 
$\{ z_{i,j} \}$ $( 1\leq i < j \leq n)$ are defining equations 
for the toric variety  $P(2,n)$ in ${\bf P}^{{ n \choose 2}-1}$:  
\[  z_{i_1,i_4}z_{i_2,i_3} -  z_{i_1,i_3}z_{i_2,i_4} = 0,  \;\; 
( 1 \leq i_1 < i_2 < i_3 < i_4 \leq n). \]
} 
\end{example}

The following theorem is due to B. Sturmfels (\cite{St}, Prop. 11.10.)

\begin{theorem} 
There exists a natural flat deformation of the 
Pl\"ucker-embedded Grassmannian 
$$G(k,n) \subset {\bf P}^{{ n \choose k}-1}$$
whose special fiber is isomorphic to the subvariety defined 
quadratic homogeneous binomial equations
\[  z_a z_{a'} - z_{min(a,a')}z_{max(a,a')}  \]
for all pairs $(a,a')$ of noncomparable elements $a, a' \in 
A(k,n)$. 
\end{theorem}   

\begin{corollary} 
The toric variety $P(k,n) \subset {\bf P}^{{ n \choose k}-1}$ 
is isomorphic to a flat degeneration of the Pl\"ucker 
embedding of the Grassmannian $G(k,n)$. 
\label{def-gr}
\end{corollary}

\section{Equations for Mirror Manifolds} 

\subsection{The mirror construction} 

Recall the definition of nef-partions for Gorenstein toric 
Fano varieties and the mirror construction 
for Calabi-Yau complete intersections  associated with 
nef-partitions \cite{LB} (we will follow the notations in \cite{BB2}).  

\begin{definition}
{\rm Let $\Delta \subset M_{\bf R}$ be a reflexive polyhedron, 
$\Delta^* \subset N_{\bf R}$ its  
dual, and $\{ e_1, \ldots ,e_l\}$ the set of vertices 
of $\Delta^*$ corresponding to torus invariant divisors 
$D_1, \ldots, D_l$ on the Gorenstein toric Fano 
variety ${\bf P}_{\Delta}$. We set $I := \{ 1, \ldots, l \}$.  
A  partition  $I= J_1 \cup \cdots \cup J_r$  of $I$ into 
a disjoint union of subsets $J_i \subset I$ 
is called a {\bf nef-partition}, if 
\[  \sum_{j \in J_i} D_j  \]
is a semi-ample Cartier divisor on ${\bf P}_{\Delta}$ 
for all $i =1, \ldots, r$.}
\end{definition}

\begin{definition} 
{\rm Let  $I= J_1 \cup \cdots \cup J_r$ be a nef-partition. 
We define  the polyhedron $\nabla_i$ $(i =1, \ldots, r)$ 
as the convex hull of $0 \in \Delta$ and all vertices 
$e_j$ with $j \in J_i$. By $\Delta_i \subset M_{\bf R}$  
$(i =1, \ldots, r)$ we denote the 
supporting polyhedron for global sections of the corresponding 
semi-ample invertible sheaf ${\cal O}( \sum_{j \in J_i} D_j)$ on 
${\bf P}_{\Delta}$.
For each  $i =1, \ldots, r$,  we denote by $g_i$ (resp. 
by $h_i$) a generic Laurent polynomial with the Newton polyhedron 
$\Delta_i$ (resp. $\nabla_i$). } 
\end{definition}

The mirror construction in \cite{LB} says that  the mirror 
of a compactified generic Calabi-Yau complete intersection
$g_1 = \cdots = g_r = 0$ is  a compactified generic Calabi-Yau 
complete intersection defined by the equations 
$h_1 = \cdots = h_r = 0$. 

Now we specialize the above mirror construction for the case 
$\Delta = \Delta^*(k,n)$, $\Delta^* = \Delta(k,n)$, and ${\bf P}_{\Delta} = 
P(k,n)$, where 
$\Delta(k,n)$ is a reflexive polyhedron defined in \ref{polyh},  
$\Delta^*(k,n)$ its polar-dual reflexive polyhedron and $P(k,n)$ the 
Gorenstein toric Fano degeneration of the Grassmannian $G(k,n)$.

\begin{definition}
{\rm Define the following $n$ subsets 
$E_1, \ldots, E_n$ of the set of vertices $\{ u_{i,j}, v_{i',j'} \}$ 
of the polyhedron $\Delta(k,n)$: 
\[ E_1 : =\{u_{1,0}\}, \;
 E_i = \{ u_{i,0}, u_{i,1}, \ldots, u_{i,n-k-1} \}, \; i =2, \ldots,k,  \]
\[ E_{k+j}: = \{ v_{1,j}, v_{2,j}, \ldots, 
v_{k,j} \}, \; j =1, \ldots, n-k-1, \; 
 E_n: = \{ v_{k,n-k}\}. \]} 
\end{definition}

\begin{proposition} 
Let $D(E_i) \subset P(k,n)$ $(i =1, \ldots, n)$ be the torus 
invariant divisor  whose irreducible components have multiplicity $1$ 
and correspond to vertices of $\Delta(k,n)$ from the subset $E_i$. 
Then the class of $D(E_i)$ is an ample generator of $Pic(P(k,n))$. 
\label{gener}
\end{proposition} 

\noindent
{\em Proof.}   By a direct computation, one obtains that 
for all $i, j \in \{ 1, \ldots, n \}$  the difference 
$D(E_i) - D(E_j)$ is a principal divisor, i.e, 
all divisors $D(E_1), \ldots, D(E_n)$ are linearly equivalent. 
On the other hand, 
\[ D(E_1) +  \cdots +  D(E_n) \]
is the ample anticanonical divisor on $P(k,n)$. By \ref{gl-sections}, 
the anticanonical divisor on $P(k,n)$ is linearly equivalent 
to $nH$, where $H$ is an  ample generator of $Pic(P(k,n))$. 
Hence, each divisor $D(E_i)$ is linearly equivalent to $H$. 
\qed

\begin{definition} 
{\rm Let $1 \leq d_1 \leq \cdots \leq d_r$ be positive integers 
satisfying the equation 
$$d_1 + \cdots + d_r = n$$ and 
 $I:= \{ 1, \ldots, n \}$. We denote by 
$X:= X_{d_1, \ldots, d_r} \subset G(k,n)$  a Calabi-Yau 
complete intersection of hypersurfaces of 
degrees $d_1, \ldots, d_r$ with $G(k,n) \subset 
{\bf P}^{ { n \choose k } -1}$.   Consider a partition 
 $I= J_1 \cup \cdots \cup J_r$  of $I$ into 
a disjoint union of subsets $J_i \subset I$ with $|J_i| = d_i$. } 
\label{J's}
\end{definition}

\begin{definition}
{\rm Let $\nabla_{J_i}$  $(i =1, \ldots, r)$ be the 
convex hull of $0 \in N(k,n)_{\bf R}$ and all vertices 
of $\Delta(k,n)$ contained in the union 
\[ \bigcup_{j \in J_i } E_j. \]
We denote by $h_{J_i}(X)$ a generic Laurent polynomial 
in variables  $X_{i',j'} := X^{f_{i',j'}}$ $( 1 \leq i' \leq k, \; 
1 \leq j' \leq n-k)$ having  $\nabla_{J_i}$ as a   Newton polyhedron.  }
\label{nablas}
\end{definition}

By \ref{gener}, one immediately obtains the following: 

\begin{corollary} 
Let $Y:=  Y_{d_1, \ldots, d_r} \subset P(k,n)$  a Calabi-Yau 
complete intersection of hypersurfaces of 
degrees $d_1, \ldots, d_r$ with the Gorenstein toric Fano 
variety $P(k,n) \subset {\bf P}^{ { n \choose k } -1}$. Then the mirror 
$Y^*$ of $Y$ $($according to  \cite{BS} and  \cite{LB}$)$ 
is a  compactified generic Calabi-Yau complete intersection 
defined by the equations 
\[ h_{J_1}(X) = \cdots =  h_{J_r}(X) = 0. \]
\end{corollary}

\begin{definition} {\rm Define $n$ 
Laurent polynomials in $k \times (n-k)$ variables $X_{i,j} := X^{f_{i,j}}$ 
as follows: 
\[ p_1(X) = a_{1,0} X^{u_{1,0}}, \;
 p_i(X) = \sum_{j =0}^{n-k-1} a_{i,j}X^{u_{i,j}}, \; i =2, \ldots,k,  \]
\[ p_{k+j}(X) = \sum_{i =1}^{k} b_{i,j}X^{v_{i,j}}, \; j =1, \ldots, n-k-1, \; 
 p_n(X) = b_{k,n-k}X^{v_{k,n-k}}, \]
where $a_{i,j}$ and $b_{l,m}$ are generically choosen complex numbers. 
In particular, the Newton polyhedron of $p_i(X)$ is the convex hull of 
$E_i$. } 
\end{definition}

\begin{conjecture} 
Let $I = \{1, \ldots, n\} =  J_1 \cup \cdots \cup J_r$ be a partition 
 of $I$ into 
a disjoint union of subsets $J_i \subset I$ with $|J_i| = d_i$ as in 
$($\ref{J's}$)$ and $Y^*_0$ be a Calabi-Yau compactification of 
a general  complete intersection 
in $({\bf C}^*)^{k(n-k)}$ defined by the equations
\[ 1 - \sum_{j \in J_i} p_j(X) = 0 \;\; ( i =1, \ldots, n),  \]  
where the coefficients 
$a_{i,j}$ and $b_{l,m}$ satisfy the 
following $(k-1)(n-k-1)$ conditions
\[ a_{k+1 -i, j-1} b_{k+1 -i,j} = a_{k+1 -i,j}b_{k-i,j}.  \]
Then a minimal desingularization $X^*$ of $Y_0^*$ 
is a mirror of a generic Calabi-Yau complete 
intersection $X:= X_{d_1, \ldots, d_r} \subset G(k,n)$. 
\label{mirror-c}
\end{conjecture} 

\begin{example} 
{\rm If $X: =X_{1,1,3} \subset G(2,5)$, we take $J_1 = \{1\}$, 
$J_2 = \{5\}$ and $J_3 = \{2,3,4 \}$. Then the mirror construction for $X$ 
proposed by \ref{mirror-c} coincides with the one considered 
in 2.2. } 
\end{example}

\subsection{Lax operators of Grassmannians}

In the paper \cite{EHX} Eguchi, Hori, and Xiong have computed 
the Lax operator $L$ for various Fano manifolds V: 
projective spaces, Del Pezzo surfaces and Grassmannians. In particular
for $V = {\bf P}^n$ the corresponding Lax operator $L$ is given by the 
formula:
\[ L = X_1 + X_2 + \cdots + X_n + qX_1^{-1} X_2^{-1} 
\cdots X_n^{-1},  \]
where $\log q $ is an element of $H_2({\bf P}^n)$. 
On the other hand, if $Z$ is an affine hypersurface 
defined by the equation 
$ L(X_1, \ldots, X_n) = 1$ 
in the algebraic torus 
$ T \cong  ({\bf C}^*)^n = {\rm Spec}\, 
{\bf C} [ X_1^{\pm 1}, \ldots, X_n^{\pm 1}]$,
then, according to \cite{Ba2}, 
a suitable compactification of $Z$ is a Calabi-Yau variety which 
is mirror dual to Calabi-Yau hypersurfaces of degree $n+1$ in ${\bf P}^n$. 

\begin{remark} 
{\rm It is natural to suggest that the last observation
can be used  
as a guiding principle for the construction of mirror manifolds 
of Calabi-Yau hypersurfaces $X$ in Fano manifolds $V$.}
\end{remark} 

Let $V$ be a Fano manifold of dimension $n$. Denote by $P$ (resp. by 
$[V]$) the class of 
unity (resp. the class of the normalized by unity volume form on $V$) 
in the cohomology ring $H^*(V)$. Let 
\[ \omega =  \frac{dX_1}{X_1} \wedge \cdots \wedge \frac{dX_n}{X_n} \]
be the invariant differential $n$-form on the $n$-dimensional algebraic 
torus $T \cong  ({\bf C}^*)^n$. According to 
\cite{EHX},   the Lax operator $L(X)$ of the Fano manifold $V$
is a Laurent polynomial in $X_1, \ldots, X_n$ with coefficients 
in the group algebra ${\bf Q}[ H_2(V,{\bf Z})]$ satisfying 
for all $m \geq 0$ the equation  
\[ \langle \sigma_m([V])P) \rangle = \frac{1}{m+1} 
\int_{\gamma} L^{m+1}(X) \omega. \]
where $\sigma_m([V])$ is the $m$-gravitational descendent of $[V]$ on 
the moduli spaces of stable maps of curves of genus $g =0$ to $V$, 
$ \langle \sigma_m([V])P) \rangle$ is the corresponding two point 
correlator function, and $\gamma$ is the standard 
generator of $H_n(T, {\bf Z})$. 

For the case $V= G(r,s)$ $(n = r(s-r))$ the following was conjectured
in  \cite{EHX}: 

\begin{conjecture} 
The Lax operator of the Grassmannian $G(r,s)$ has the following form
\[ L(X) = X_{[1,1]} + 
\sum_{\begin{array}{c} {\scriptstyle  1 \leq a \leq s-r
 } \\
{\scriptstyle  1 \leq b \leq r } \end{array} } 
X_{[a,b]}^{-1}(X_{[a+1,b]} + X_{[a,b+1]}) 
+ q X_{[s-r,r]}^{-1}, \]
where $\log q  \in H_2(G(r,s))$ and $X_{a,b} = 0$ if $a > s-r$ or $b > r$.
\label{lax}  
\end{conjecture}

\begin{proposition}
Let $P(r,s)$ be the toric degeneration of the 
Grassmannian $G(r,s)$. Then the  equation $L(X) = 1$ defines  
a $1$-parameter subfamily in the family of 
toric mirrors of Calabi-Yau hypersurfaces in $P(r,s)$ 
$($see {\rm \cite{Ba2}}$)$.
\end{proposition} 

\prf According to \cite{Ba2}, we have to identify the Newton 
polyhedron of the Laurent polynomial $L(X)$  in Conjecture \ref{lax} with 
the reflexive polyhedron $\Delta(r,s)$. The latter follows immediately 
from the explicit description of $\Delta(r,s)$ in \ref{polyh} and from the 
$1$-to-$1$-correspondence 
$f_{i,j} \leftrightarrow X_{[j,i]}$.\qed 
\vskip 10pt

\begin{proposition} 
The equations for the mirrors to  
Calabi-Yau hypersurfaces conjectured in \ref{mirror-c} in 
$G(r,s)$ coincide with the equations $L(X) =1$ where 
$L(X)$ is the Lax operator conjectured for $G(r,s)$ 
in \cite{EHX}.  
\end{proposition}

\noindent
\prf  It is easy to see that the coefficients of the polynomial 
$L(X)$ satisfy all $r(s-r)$ monomial relations which reduce to the 
equality $1 \cdot 1 = 1 \cdot 1$. On the other hand, using 
the action of the $r(s-r)$-dimensional torus  on the coefficients of the 
Laurent polynomial 
\[ 1 - (p_1(X) + \cdots + p_s(X))  \]
defining the mirror in  \ref{mirror-c}, one can reduce to only one independent 
parameter, for instance, the unique coefficient $b_{r,s-r}$ of $p_s(X) = 
b_{r,s-r}X_{r,s-r}^{-1}$. By setting $q : = b_{r,s-r}$ and 
$X_{[i,j]}: = X_{j,i}$, we can identify the variety $Y_0^*$ in 
\ref{mirror-c} with a toric compactification of the affine hypersurface 
$L(X) =1$. 
\qed 
\vskip 10pt

Using the explicit description  of the multiplicative 
structure of the small quantum cohomology of $G(k,n)$, 
it is not difficult to check Conjecture \ref{lax} for each 
given $r$ and $s$:

\begin{example} 
{\rm The Lax operator of the  
Grassmannian $G(2,4)$ is 
\[  X_{[1,1]} + X_{[1,1]}^{-1}(X_{[2,1]} + X_{[1,2]}) + 
X_{[2,1]}^{-1}X_{[2,2]} + X_{[1,2]}^{-1}X_{[2,2]}
 + q X_{[2,2]}^{-1}. \]
Its Newton polyhedron is isomorphic to $\Delta(2,4)$ from 2.1.} 
\end{example}

\begin{example} 
{\rm  The Lax operator of the  
Grassmannian $G(2,5)$ is 
\[  X_{[1,1]} + 
X_{[1,1]}^{-1}(X_{[2,1]} + X_{[1,2]})  + 
X_{[2,1]}^{-1}( X_{[3,1]} + X_{[2,2]} ) + 
X_{[1,2]}^{-1} X_{[2,2]} + 
X_{[2,2]}^{-1} X_{[3,2]} + q X_{[3,2]}^{-1}. \]
Its Newton polyhedron is isomorphic to $\Delta(2,5)$ from 2.2.} 
\end{example}

\section{Hypergeometric series} 

\subsection{The Trick with the Factorials}\label{trick} 

If $X$ is a Calabi-Yau the complete intersection of 
hypersurfaces of degree $l_1, l_2, \ldots, l_r$ in ${\bf P}^n$, then the 
generalised hypergeometric series
$$ \Phi_X(q)=\sum_{m=0}^{\infty} 
\frac{(l_1m)!(l_2m)!\ldots(l_rm)!}{(m!)^{n+1}}q^m$$ 
is  main period of its  mirror $X^*$. 
As is well-known, one can obtain the instanton numbers 
for $X$ by a formal manipulation with this series,
 see e.g. \cite{BS} and \ref{instanton}. More precisely, one 
transforms the Picard-Fuchs differential operator $P$ annihilating the series
$\Phi_X$ to the form $D^2 \frac{1}{K(q)}D^2$
(where $D=q\partial/\partial q$) and reads off 
the the $n_d$ from the power series expansion of the function $K$:\\
$$K(q)=l_1 l_2 \ldots l_r +\sum_{d=1}^{\infty} n_d d^3 \frac{q^d}{1-q^d}.$$

It is important to observe that the power 
series  $\Phi_X$ can be obtained from a power series 
$$A_V(q)=\sum_{m=0}^{\infty} \frac{1}{(m!)^{n+1}}q^m$$
by the multiplication of its $m$-th 
coefficient by the product 
of factorials $(l_1m)!(l_2m)!\ldots(l_rm)!$.
On the other hand, the power series $A_V$ can be characterized 
as the unique series  
$A_V=1+\ldots$ solving the differential equation
$((q\frac{\partial}{\partial q})^{n+1}-q)A_V=0$
associated with the small quantum cohomology of ${\bf P}^n$. 
This differential equation arizes as the reduction of the first order
differential {\em system}
$$q\frac{\partial}{\partial q} \vec{S} = p \circ \vec{S}$$
for a $H^*({\bf P}^n)$-valued function $\vec{S}=S_0+S_1p+\ldots+S_np^n,$
where  $p \in H^2({\bf P}^{n+1})$ is an ample  generator, 
$\{1, p, p^2, \ldots, p^n\}$ 
is  a basis for $H^*({\bf P}^n)$, and  $p \circ$ is the operation of
{\em quantum multiplication} with $p$ in the small quantum cohomology of 
${\bf P}^{n}$. Since it is well-known that the small 
quantum cohomology ring of ${\bf P}^n$ is 
defined by the relation $(p \circ)^{n+1} -q=0$, one finds immediately 
comes to the differential equation.
In particular, we see that the function $A_V$ is uniquely determined 
by the small quantum cohomology ring of $V={\bf P}^n$.

It is natural to try to use these ideas to obtain $\Phi_X$ from $A_V$ 
for varieties other than ${\bf P}^n$, for example 
for Grassmannians or other Fano varieties. 
If it works, this method allows one to 
find instanton numbers without knowing an explicit mirror manifold. 
We will formulate this trick in
some generality below.

Let $V$ be a smooth projective variety, which for reasons of
simplicity of exposition is assumed to have only even cohomology and
and that $H^2(V,\integer ) \cong H_2(V,\integer) \cong {\bf Z}$. Let
$p$ be the ample generator of $H^2(V, \integer)$, 
$\gamma$ a positive generator
for $H_2(V, \integer)$. We denote by $1_V \in 
H^0(V)$ the fundamental class of $V$ and by $<-,->$ 
the Poincar\'e pairing. The small quantum cohomology ring
$QH^*(V)$ of $V$ is the free $\rational[[q]]$-module $H^*(V,\rational[[q]])$
with a new multiplication $\circ$ determined by
$<A \circ B,C>=<A,B,C>=\sum_{m=0}^{\infty}<A,B,C>_mq^m>$
where $$<A,B,C>_m=I^V_{0,3,m\gamma}=\int_{[\overline{M}_{0,3}]}e^*_1(A) 
\cup e_2^*(B) \cup e^*_3(C)$$
are the {\em $3$-point, genus 0, Gromov-Witten 
invariants}, see  \cite{FP}. The operator of quantum multiplication 
with the ample generator $p \in H^2(V, \integer)$ defines the 
{\em Quantum Differential System}, see e.g. \cite{G1}:
$$\frac{\partial}{\partial t}\vec{S}=p \circ \vec{S}$$
where $\vec{S}$ is an series in the variable 
$t=log $ with coefficients from $H^*(V,\rational)$. 
The {\em Quantum Cohomology $\cD$-module}is the $\cD$-module 
generated by the top components 
$<\vec{S}, 1_V>$ of all solutions $\vec{S}$ to the above differential system.
In the case under consideration, it will be of the form $\cD/\cD P$, for a
certain differential operator $P$.

\begin{definition}{\em
The {\em $A$-series of $V$} is the unique solution of the 
Quantum Cohomology  $\cD$-module of the form
$A_V=\sum_{m=1}^{\infty} a_m q^m$
with $a_0=1$.
 }
\end{definition} 

Let $X$ be the intersection of hypersurfaces  
of degree $l_1, l_2, \ldots, l_r$ in $V$. 
In other words, $X$ is the zero-set
of a generic section of the decomposable bundle 
$\cE:=\cO(l_1p)\oplus \cO(l_2p)\oplus\ldots\oplus\cO(l_rp)$.

\begin{definition}
{\em Let  $A_V=\sum_{m=1}^{\infty} a_m q^m$ be the 
$A$-series of a Fano manifolds $V$. 
Define the {\em ${\cal E}$-modification of $A_V$} as 
follows:
$$\Phi_{\cE}(q): =\sum_{m=0}^{\infty} a_m \prod_{i=1}^r (m l_i)! q^m.$$
}
\end{definition}

\begin{definition}
{\em Assume that $X \subset V$ has trivial canonical class, i.e, 
$X$ is a Calabi-Yau variety.  
We say that the {\em Trick with the Factorials works}, if the 
function $\Phi_{\cE}$ is equal to the monodromy invariant period 
$\Phi_X$ of the mirror family $X^*$ in some algebraic parametrization.}
\end{definition}

If the Trick with the Factorials works,  
the usual formal manipulation (see \cite{BS}, \ref{instanton}) with
the series $\Phi_{\cE}$ produces the instanton numbers for $X$!\\

\begin{remark}
{\rm  (i)  It is possible to formulate the Trick with the Factorials in much
greater generality  \cite{K3}, \cite{BCKS}.

 (ii)  The $A$-series $A_V$ very well can be identically $1$, 
but if $V$ is Fano, it will contain interesting information 
and it is in such cases that the Trick with the Factorials 
has a chance to work.

 (iii)  A better formulation uses instead of $A_V$ a certain
cohomology-valued series $S_V$, whose components  make up a complete
solution set to the quantum $\cD$-module. Instead of the
factorially modified series $\Phi_{\cE}$ 
one has a factorially modified cohomological function $F_{\cE}$. 
We say that Trick with the Factorials works, if $S_V$ and
$F_{\cE}$ differ by a coordinate change  \cite{K3}, \cite{BCKS}.
Such a theorem is a form of the  Lefschets hyperplane section theorem
in quantum cohomology.
 
 (iv)  Givental's mirror theorem for toric varieties, \cite{G3}, 
implies that the
Trick with the Factorials works for complete intersections in toric
varieties. 

 (v)  More generally, it follows from a recent theorem of Kim, 
\cite{K3}, that the
Trick with the Factorials works for arbitrary homogeneous spaces.

(vi)  E. Tj{\o}tta has applied the Trick with the Factorials 
succesfully in a
non-homogeneous case, \cite{tjotta}.}
\end{remark}

\subsection{Hypergeometric solutions for Grassmannians}

In this paragraph we apply the above ideas to the case of
Grassmannians.

In \cite{BCKS}, we describe a simple rule for writing down the 
$GKZ$-hypergeometric series $A_{P(k,n)}$ assiciated with the 
Gorenstein toric Fano variety 
$P(k,n)$  in terms of the combinatorics of a certain 
graph. Here we give a formula for $A_{P(k,n)}$ without going into 
the details:

$$A_{P(k,n)}(q,\tilde{q})=
\sum_{s_{i,j}\ge 0}\frac{1}{(m!)^{n}}
\prod_{i=1}^{k-1}\prod_{j=1}^{n-k-1}{s_{i+1,j}
\choose s_{i,j}}{s_{i,j+1}\choose s_{i,j}}q^m \tilde{q}_{i,j}^{s_{i,j}}$$
where we put $s_{i,j}=m$ if $i > k-1$ or $j > n-k-1$.

\begin{example} 
{\rm $G(2,5)$: 
$$ A_{P(2,5)}(q, \tilde{q})=
\sum_{m,r,s \ge 0}\frac{1}{(m!)^5}{m \choose r}{s \choose r}
{m \choose s}^2q^m \tilde{q}_1^r\tilde{q}_2^s.$$
} 
\end{example}
\vskip 10pt

\begin{example}
{\rm 
$G(3,6)$ : 
$$ A_{P(3,6)}(q, \tilde{q})=
\sum_{m,r,s,t,u}\frac{1}{(m!)^6}{s \choose r}{t \choose r} 
{m \choose s}{u \choose s}{u \choose t}{m \choose t}
{m \choose u}^2 q^m \tilde{q}_1^r \tilde{q}_2^s \tilde{q}_3^u \tilde{q}_4^v.$$
}
\end{example}
\vskip 10pt

We conjecture an  explicit general formula 
for the  series $A_{G(k,n)}(q)$ of an arbitrary Grassmannian: 

\begin{conjecture} 
Let $A_{P(k,n)}(q, \tilde{q})$ be the 
A-hypergeometric series of the toric variety $\widehat{P(k,n)}$ as above. Then 
\[A_{G(k,n)}(q) = A_{P(k,n)}(q, {\bf 1}). \]  
\label{flagmirror} 
\end{conjecture}

Using the explicit formulas for multiplication in the quantum cohomology
of Grassmannians  \cite{Ber}, one can write down the 
Quantum Differential System for $G(k,n)$ and reduce this 
first order system to a higher 
order differential equation satisfied by its components. In particular, one
can write down the differential operator $P$ annihilating the
component $<\vec{S},1>$ of any solution $\vec{S}$.

Below we record some of the (computer aided) 
calculations of the operator $P$ we did ($D$ denotes the 
operator $\partial/\partial t=q \partial/\partial q$):

$$
\begin{array}{|lcl|}
\hline
 & & \\
G(2,4)&:&D^5-2q(2D+1)\\ 
 & & \\
G(2,5)&:&D^7(D-1)^3-qD^3(11D^2+11D+3)-q^2 \\ 
 & & \\
G(2,6)&:&D^9(D-1)^5-qD^5(2D+1)(13D^2+13D+4)\\ 
 & & \\
 & &-3q^2(3D+4)(3D+2)\\ 
 & & \\
G(3,6)&:&D^{10}(D-1)^4-qD^4(65D^4+130D^3+105D^2+40D+6)\\ 
 & & \\
 & &+4q^2(4D+3)(4D+5)\\
 & & \\
\hline
\end{array} 
$$
The operator for $G(2,7)$ is:
$$
\begin{array}{c}
D^{11}(D-1)^7(D-2)^7(D-3)^7(D-4)^3\\ \\
-\frac{1}{3}qD^7(D-1)^7(D-2)^7(D-3)^3(173D^4+340D^3+272D^2+102D+15)\\ \\
-\frac{2}{9}q^2D^7(D-1)^7(D-2)^3(1129D^4+5032D^3+7597D^2+4773D+1083)\\ \\
+\frac{2}{9}q^3D^7(D-1)^3(843D^4+2628D^3+2353D^2+675D+6)\\ \\
-\frac{1}{9}q^4D^3(295D^4+608D^3+478D^2+174D+26)+\frac{1}{9}q^5,\\
\end{array}
$$
while the one for $G(2,8)$ takes about two pages. 
Clearly, since both the structure of the quantum cohomology 
ring and the hypergeometric series are very explicit, one 
should seek a better way to prove Conjecture \ref{flagmirror}.

Nevertheless, using the above operators one obtains by direct computation 
the following: 

\begin{theorem} The conjecture \ref{flagmirror} is true for
$G(2,4)$, $G(2,5)$, $G(2,6)$, $G(2,7)$, $G(3,6)$. \qed 
\label{grasscheck}
\end{theorem}

\section{Complete Intersection  Calabi-Yau $3$-folds}

\subsection{Conifold transitions and mirrors}

Now we turn our attention to the main point of the 
paper, namely the construction, via 
conifold transitions, of mirrors for Calabi-Yau $3$-folds $X$ 
which are complete 
intersections in Grassmannians $G(k,n)$.

By \ref{sing-l}, the singular locus of a generic $3$-dimensional complete intersection $X_0$ 
of $P(k,n)$ with 
$r$ hypersurfaces $H_1, \ldots, H_r$ of degrees $d_1, \ldots, d_r$ 
($[H_i] = d_i[H]$, $i =1, \ldots, r$) consists of 
\[ p = d_1 d_2 \cdots d_r \left( \sum_{i =1,\, j =1}^{k-1,\,n-k-1} 
d(W_{i,j}) \right) \]
nodes, where $d(W_{i,j})$ is the degree of $W_{i,j}$ with respect to 
the generator $H$ of the Picard group of $P(k,n)$.    
On the other hand, by \ref{def-gr}, $X_0$ is a flat degeneration 
of the smooth Calabi-Yau $3$-fold $X \subset G(k,n)$. 

The small crepant resolution $\widehat{P}(k,n) \lra P(k,n)$ of 
the ambient toric variety
induces a small crepant resolution $Y\lra X_0$. Hence $Y$ is a {\em smooth}
Calabi-Yau complete intersection in the toric variety $\widehat{P}(k,n)$, 
which is obtained from $X$ by a {\em conifold transition}.

\begin{theorem} Let $p$ be the number of nodes 
of $X_0$, and let $\alpha = (k-1)(n-k-1)$. 
Then the Hodge numbers of $X$ and $Y$ are related by
$$h^{1,1}(Y)=h^{1,1}(X)+\alpha$$
and $$h^{2,1}(Y)=h^{2,1}(X)+\alpha -p.$$
\label{formul-2}
\end{theorem}

\noindent
\prf By construction, $Y$ is a complete intersection of general sections of 
big {\it semiample} line bundles on $\widehat{P}(k,n)$ 
(i.e., line bundles which are generated by global 
sections and big). Using the explicit formula for 
$h^{1,1}(Y)$ from (\cite{BB}, 
Corollary 8.3) and the fact that the only boundary lattice points 
of $\Delta(k,n)$ are its vertices, we obtain the isomorphism 
 $\Pic(Y)\cong\Pic(\widehat{P})$, 
which gives the first relation. 
On the other hand, the $p$ vanishing 
3-cycles on X that shrink to nodes in 
the degeneration must satisfy $\alpha$ linearly 
independent relations by \cite{C}, and the second relation follows. \qed

\vskip 10pt

The mirror construction for complete intersection 
Calabi-Yau manifolds in toric
varieties given in \cite{Ba2, BB} provides us 
with the mirror family of Calabi-Yau manifolds $Y^*$. 
The generic member of this family is nonsingular
(it is obtained by a MPCP-resolution of the ambient toric variety).
There is a natural isomorphism of the Hodge groups 
$H^{1,1}(Y)\lra H^{2,1}(Y^*)$ (see \cite{Ba2, BB}). 

During the conifold transition from $X$ to $Y$, 
we have increased the "K\"ahler moduli", that is, 
the rank of $H^{1,1}$. This says that we should really 
look at the one-parameter subfamily of mirrors given by the subspace 
of $H^{2,1}(Y^*)$ corresponding via the isomorphism above 
to the divisors on $Y$ which come from $X$.
For this reason, the generalized hypergeometric series 
$\Phi_X$ of $X^*$ is a specialization
of the monodromy invariant period integral of the mirror family $Y^*$
to the subfamily $Y_0^*$ defined in \ref{mirror-c}.

Let $\nabla_{J_1}, \ldots, \nabla_{J_r}$ be convex polyhedra as in 
\ref{nablas}. Denote by $\nabla(k,n)$ the Minkowski sum of 
$\nabla_{J_1}, \ldots, \nabla_{J_r}$. Then  $\nabla(k,n)$ is 
a reflexive polyhedron and ${\bf P}_{\nabla(k,n)}$ is 
a Gorenstein toric Fano variety 
defined by a nef-partition corresponding to  the equation 
\[ d_1 + \cdots + d_r = n. \]

\begin{conjecture}
After a MPCP-desingularization of the ambient toric variety 
${\bf P}_{\nabla(k,n)}$, 
the general member $Y_0^*$ of the special $1$-parameter subfamily 
is a Calabi-Yau variety with the same number $p$ 
of nodes as $X_0$, satisfying
 $\alpha -p$ relations. A small resolution 
$X^*$ of $Y^*_0$ is a mirror of $X$. 
\label{toricmirror}
\end{conjecture}

\begin{remark} 
{\rm The statement \ref{toricmirror} has been easily  checked for the two 
simplest cases of Section 2, where the toric mirror construction reduces 
to a hypersurface case.  However, singularities of  $Y_0^*$ 
are more difficult control for $4$ remaining cases which can not be 
reduced  to Calabi-Yau hypersurfaces in $4$-dimensional Gorenstein 
toric Fano varieties.} 
\end{remark}

\subsection{The computation of instanton numbers}\label{instanton}

We denote  
by $X_{d_1,\ldots, d_r}\subset G(k,n)$ a 
Calabi-Yau complete 
intersection of $r$ hypersurfaces of degrees $d_1, \ldots, d_r$ 
with the Grassmannian $G(k,n) \subset {\bf P}^{{n \choose k} -1}$.   
We denote by $Y$ the toric Calabi-Yau complete intersection 
in $\widehat{P}(k,n)$ obtained by 
a conifold transition via resolution of $p$ nodes on the 
degeneration $X_0$ of $X$ $(h^{1,1}(X) = 1$), and by 
$\alpha$ the number of relations satisfied by 
the homology classes of the corresponding $p$ vanishing $3$-cycles on $X$.

Now we  list all  cases of Calabi-Yau complete intersection
$3$-folds $X$ in Grassmannians and collect the information about 
topological invariants of  $X$ 
and their conifold modifications $Y$. 

\begin{center}

\begin{tabular}{|c|c|c|c|c|c|c|c|} \hline
$X$ & $h^{2,1}(X)$ & 
$\chi(X)$ & $h^{1,1}(Y)$ & $h^{2,1}(Y)$ & $\chi(Y)$ & $\alpha$ & $p$ \\ 
\hline 
$X_4 \subset G(2,4)$ & $89$ & $- 176$ & $2$ & $86$ & $-168$ & $1$ & $4$ \\
\hline 
$X_{1,1,3} \subset G(2,5)$ & $76$ & $-150$ & $3$ & $72$ & $-138$ & $2$ 
& $6$ \\
\hline
$X_{1,2,2} \subset G(2,5)$ & $61$ & $-120$ & $3$ & $55$ & $ -104$ & $2$ 
& $8$ \\
\hline
$X_{1,1,1,1,2} \subset G(2,6)$ & $59$ & $-116$ & $4$ & $52$  & $-96$ & $3$ 
& $10$ \\
\hline
$X_{1, \ldots, 1} \subset G(2,7)$ & $50$ & $ -98$ & $5$ & $40$ & $-70$ & $4$ 
& $14$ \\
\hline
$X_{1, \ldots, 1} \subset G(3,6)$ & $49$ & $-96$ & $5$ & $37$ & $-64$ & $4$ & 
$16$ \\         
\hline
\end{tabular}
\end{center}

Recall the (standard) formal procedure used  to compute 
the instanton numbers. (More details can 
be found e.g in \cite{BS}.)
We set 
$$\Phi_X(z) := \sum_{m \geq 0} b_mz^m $$
to be the generalized hypergeometric series (with variable $z$) 
corresponding to the 
monodromy invariant period of the mirror $X^*$. As was explained in 
\ref{trick}, one can start with the the $A$-series for the
grassmannian, and apply the Trick with the Factorials to find the
coefficients $b_m$.

Then $\Phi_X(z)$ satisfies a Picard-Fuchs differential equation 
\[ P \Phi_X(z) = 0, \]
where 
$P$ is a differential operator of order $4$ having a maximal unipotent 
monodromy at $z= 0$. We compute $P$ by finding an explicit 
recursion relation among coefficients $b_m$ of the 
generalized hypergeometric series $\Phi_X(z)$. To bring $P$ into the
form $D^2 \frac{1}{K}D^2$, one has to 
change the coordinate $z$ to $q = {\rm exp}(\Phi_1(z)/\Phi_X(z))$,
where $\Phi_1$ is the logarithmic solution of $P$. To obtain $K$, it is
convenient to use the Yukawa coupling.  In the coordinate $z$ it has form 
\[ K_{zzz} =  
 \frac{K_z^{(3)}}{\Phi_X^2(z)} \left( \frac{dz}{z} \right)^{\otimes 3}, \]
where $K_z^{(3)}$ is some rational function of $z$ that can be determined directly 
from $P$. The Yukawa coupling in coordinate $q$ then is of the form 
 \[ K_{qqq} =  K_q^{(3)} 
\left( \frac{dq}{q} \right)^{\otimes 3}, \]
where 
\[ K_q^{(3)} = n_0 + \sum_{m =1}^{\infty} n_m \frac{m^3 q^m}{1 -q^m} \]
and $n_m$ are the instanton numbers for rational curves  of degree 
$m$ on $X$.

From proposition \ref{grasscheck} and 
Kim's  Quantum Hyperplane Theorem (\cite{K3}), we have the following

\begin{theorem}
The virtual numbers of rational curves on a general complete intersection
Calabi-Yau three-fold in a Grassmannian are the ones listed in the 
tables of the next section. 
\end{theorem}

\section{Picard-Fuchs Operators and Yukawa Couplings} 
{\tiny
\noindent
\subsection{$X_{1,1,3} \subset G(2,5)$} 

\begin{center}
\begin{tabular}{|c|c|} \hline
& \\
$ b_m $ & ${\displaystyle 
\frac{(m!)(m!)(3m)!}{(m!)^5}\sum_{r,s}{m\choose r}{s \choose r} 
{m\choose s}^2} $ \\
& \\
\hline
& \\
${P}$ & $  D^4 - 3z (3 D +2)( 3 D + 1)
(11 D^2 + 11  D + 3) $\\
& $ - 9z^2 (3 D + 5)(3 D + 2)
(3 D + 4)(3 D + 1)$ \\
& \\
\hline
& \\
$ K_z^{(3)}$ & ${\displaystyle  \frac{15}{1-11\cdot3^3z - 3^9z^2}  }$ \\
& \\
\hline
& \\
$ n_m$  & $n_1 =540 ,\; n_2 =12555 , \; n_3 =621315 ,\; n_4 =44892765 ,\;
n_5 = 3995437590  $ \\
& \\
\hline
\end{tabular}

\end{center}
}
\noindent 
{\tiny
\subsection{$X_{1,2,2} \subset G(2,5)$} 

\begin{center}

\begin{tabular}{|c|c|} \hline
& \\
$ b_m $ & ${\displaystyle \frac{ (m!)(2 m)!)^2}{(m!)^5} \sum_{r,s }
{m\choose r}{s \choose r}{m\choose s}^2}$ \\
& \\
\hline
& \\
${P}$ & $   D^4 - 4z (11  D^2 + 11  D + 3)
(1 + 2 D)^2$\\
& $ - 16z^2 (2 D + 3)^2 (1 + 2 D)^2$ \\
& \\
\hline
& \\
$ K_z^{(3)}$ & ${\displaystyle  \frac{20}{1 - 11 \cdot 2^4z - 2^8z^2}  }$ \\
& \\
\hline
& \\
$ n_m$  & $n_1 =400 ,\; n_2 =5540 , \; n_3 = 164400,\; n_4 =7059880 ,\;
n_5 = 373030720 $ \\
& \\
\hline
\end{tabular}

\end{center}
}
{\small
The locus of conifold singularities in the toric variety $P(2,5)$ consists 
of $2$ codimesion-$3$ toric strata  of degree $1$.  This gives $6$ nodes 
on the generic complete intersection of type $(1,1,3)$ in 
$P(2,5) \subset {\bf P}^{10}$ and $8$ nodes on the generic 
complete intersection of type $(1,2,2)$ in 
$P(2,5) \subset {\bf P}^{10}$    
}
\noindent
\subsection{$X_{1,1,1,1,2} \subset G(2,6)$} 
{\tiny
\begin{center}

\begin{tabular}{|c|c|} \hline
& \\
$ b_m $ &   ${\displaystyle \frac{ (m!)^4(2 m)!)}{(m!)^6} \sum_{r,s,t}
{m\choose r}{s\choose r}{m\choose s}{t\choose s}{m\choose t}^2}$\\
& \\
\hline
& \\
${P}$ & $  D^4 - 2z (4 + 13  D + 13  D^2)(1 +
2 D)^2 $\\
& $ -12 z^2 (3 D +2)( 2 D + 3)(1 +  2 D)( 3 D + 4)$ \\
& \\
\hline
& \\
$ K_z^{(3)}$ & ${\displaystyle  \frac{28}{1-26\cdot2^2z - 
27\cdot 2^4z^2}  }$ \\
& \\
\hline
& \\
$ n_m$  & $n_1 =280 ,\; n_2 =2674 , \; n_3 =48272 ,\; n_4 = 1279040,\;
n_5 = 41389992  $ \\
& \\
\hline
\end{tabular}

\end{center}
}
{\small
The locus of conifold singularities in the toric variety $P(2,6)$ consists 
of $2$ codimesion-$3$ toric strata  of degree $2$ and 
$1$ codimension-$3$ toric stratum of degree $1$. This gives $10$ nodes 
on the generic complete intersection of type $(1,1,1,1,2)$ in 
$P(2,6) \subset {\bf P}^{14}$.    
}
\noindent
\subsection{$X_{1,1,1,1,1,1,1} \subset G(2,7)$} 
{\tiny
\begin{center}

\begin{tabular}{|c|c|} \hline
& \\
$ b_m $ &   ${\displaystyle \frac{ (m!)^7}{(m!)^7} \sum_{r,s,t,u }
{m\choose r}{s\choose r}{m\choose s}{t\choose s}
{m\choose t}{u\choose t}{m\choose u}^2}$\\
& \\
\hline
& \\
${P}$ & $ 9  D^4 - 3z (15 + 102 D + 272  D^2 +
340  D^3 + 173  D^4) $\\
& $ - 2z^2 (1083  + 4773 D + 7597  D^2 + 5032  D^3 + 
1129  D^4)$ \\
& $ + 2z^3 (6 + 675  D + 2353  D^2 + 2628 D^3 + 843  D^4)$ \\
& $ - z^4(26 + 174  D + 478  D^2 + 608  D^3 + 295  D^4) + 
z^5 ( D + 1)^4 $ \\
& \\
\hline
& \\
$ K_z^{(3)}$ & ${\displaystyle  \frac{42-14z}{1-57z -289z^2 + z^3}  }$ \\
& \\
\hline
& \\
$ n_m$  & $n_1 = 196,\; n_2 =1225 , \; n_3 =12740 ,\; n_4 =198058 ,\;
n_5 = 3716944  $ \\
& \\
\hline
\end{tabular}

\end{center}
}
{\small
The locus of conifold singularities in the toric variety $P(2,7)$ consists 
of $2$ codimesion-$3$ toric strata  of degree $2$ and 
$2$ codimension-$3$ toric stratum of degree $5$. This gives $14$ nodes 
on the generic complete intersection of type $(1,1,1,1,1,1)$ in 
$P(2,7) \subset {\bf P}^{20}$.    
}

\noindent

\subsection{$X_{1,1,1,1,1,1} \subset G(3,6)$} 
{\tiny
\begin{center}

\begin{tabular}{|c|c|} \hline
& \\
$ b_m $ &  ${\displaystyle \frac{ (m!)^6}{(m!)^6} \sum_{r,s,t,u }
 {s \choose r}{t \choose r} 
{m \choose s}{u \choose s}{u \choose t}{m \choose t}
{m \choose u}^2}$\\
& \\
\hline
& \\
${P}$ & $  D^4 - z (6 + 40 D + 105  D^2 +
130  D^3 + 65  D^4) $\\
& $ +4z^2 (4 D + 5)(4 D +3)( D + 1)^2$ \\
& \\
\hline
& \\
$ K_z^{(3)}$ & ${\displaystyle  \frac{42}{1 - 65 z - 64z^2}  }$ \\
& \\
\hline
& \\
$ n_m$  & $n_1 =210 ,\; n_2 =1176 , \; n_3 =13104 ,\; n_4 = 201936 ,\;
n_5 =3824016 $ \\
& \\
\hline
\end{tabular}

\end{center}
}
{\small
The locus of conifold singularities in the toric variety $P(3,6)$ consists 
of $2$ codimesion-$3$ toric strata  of degree $2$ and 
$2$ codimension-$3$ toric strata of degree $6$. This gives $16$ nodes 
on the generic complete intersection of type $(1,1,1,1,1,1)$ in 
$P(3,6) \subset {\bf P}^{19}$.    
}

\vskip .2truein 
\noindent
\section{Acknowledgement}

\noindent
We would like to thank 
S. Katz, S.-A. Str{\o}mme, E. R{\o}dland and E. Tj{\o}tta 
for helpful discussions and 
the Mittag-Leffler Institute for hospitality. 
The second and third named authors
have been supported by Mittag-Leffler Institute postdoctoral fellowships.


\newpage

\vskip .2truein
\begin{thebibliography}{2}
         

\bibitem{ACJM} A. C. Avram, P. Candelas, D. Jancic, M. Mandelberg, 
{\em On the Connectedness of the Moduli Space of
Calabi--Yau Manifolds},  Nucl.Phys. B465 (1996) 458-472 (hep-th/9511230) 


\bibitem{AKMS} 
A.C.Avram, M.Kreuzer, M.Mandelberg and H.Skarke, {\em  The web of 
Calabi-Yau hypersurfaces in toric varieties}, hep-th/9703003  

\bibitem{Ba1} V.V. Batyrev, {\em On classifications of smooth 
projective toric varieties}, T\^ohoku  Math. J. {\bf 43} (1991), 569-585.  


\bibitem{Ba2} V.V. Batyrev, {\em Dual polyhedra and mirror symmetry 
for Calabi-Yau hypersurfaces in toric varieties}, 
J. Alg. Geom., {\bf 3} (1994) 493-535. 

\bibitem{BB} V.V. Batyrev, L. A. Borisov, {\em On Calabi-Yau Complete
Intersections in Toric Varieties},  in \lq \lq Higher
 Dimensional Complex Geometry", M. 
Andreatta and T. Peternell eds., 1996, 37-65. (alg-geom/9412017)

\bibitem{BB2} V.V. Batyrev, L.A. Borisov, {\em Dual Cones and Mirror 
Symmetry for Generalized Calabi-Yau Manifolds}, in {\sl Mirror Symmetry II}, 
(B. Greene, S-T. Yau eds.) AMS and Int. Press (1997), 71-86. 

\bibitem{BS} V. V. Batyrev, D. van Straten, 
{\em Generalized Hypergeometric Functions 
and Rational Curves on Calabi-Yau Complete Intersections in Toric 
Varieties}, Commun. Math.  Phys. {\bf 168} (1995), 493-533. 
(alg-geom/9307010) 

\bibitem{BCKS} V.V. Batyrev, I. Ciocan-Fontanine, 
B. Kim, and  D. van Straten, {\em Mirror Symmetry and Toric Degenerations 
of Partial Flag Manifolds}, Preprint 1997.  

\bibitem{Ber} A. Bertram, {\em Quantum Schubert calculus}, to appear in Adv. Math. 


\bibitem{behrend} K. Behrend, {\em Gromov-Witten invariants 
in algebraic geometry}, Invent. 
Math., {\bf 127}(3), (1997), 601-617. 


\bibitem{LB} L. A. Borisov, {\em Towards Mirror Symmetry of 
Calabi-Yau Complete Intersections in Gorenstein Toric Fano 
Varieties.} alg-geom/9310001. 

\bibitem{BLS} I. Brunner, M. Lynker, R. Schimmrigk, {\em 
 Dualities and phase transitions for Calabi-Yau 3folds and 4folds}, 
hep-th/9703182 


\bibitem{CGGK}  T.-M. Chiang, B.R. Greene, M. Gross and Y. Kanter, 
   {\em Black hole condensation and the web of Calabi - Yau 
manifolds}, hep-th/9511204 


\bibitem{CDLS} P. Candelas, A. Dale, C.A. L\"utken, R. Schimmrigk,  
{\em Complete Intersection Calabi-Yau Manifolds}, Nucl. Phys. 
{\bf B298} (1988), 493

\bibitem{CGH} P. Candelas, P. Green and T. H\"ubsch, {\em Rolling 
Among Calabi-Yau Vacua}, Nucl. Phys. {\bf B 330} (1990) 49

\bibitem{C} H. Clemens, {\em Double solids}, Adv. Math. {\bf 47} (1983), 
107-230.  

\bibitem{EHX} T. Eguchi, K. Hori and C.-S. Xiong, {\em Gravitational Quatum 
Cohomology}, hep-th/9605225 to appear in Int. J. Mod. Phys. A. 

\bibitem{FP} W. Fulton, R. Pandharipande, {\em Notes on
 Stable Maps And Quantum Cohomology}, preprint, alg-geom/9608011.

\bibitem{G1} A. Givental, {\em Equivariant Gromov-Witten Invariants}, 
IMRN, No. 13 (1996), 613-663.  
(alg-geom/9603021) 

\bibitem{G2} A. Givental, {\em Stationary Phase Integrals, Quantum 
Toda Lattices, Flag Manifolds and the Mirror Conjecture}, alg-geom/9612001. 

\bibitem{G3} A. Givental, {\em A Mirror Theorem for Toric 
Complete Intersections}, alg-geom/9701016. 

\bibitem{GL0}  V. Lakshmibai, {\em Degeneration of flag varieties to 
toric varieties}, C.R. Acad. Sci. Paris {\bf 321} (1995), 1229-1234.

\bibitem{GL1} N. Gonciulea and V. Lakshmibai, {\em Degenerations of 
Flag and
Schubert Varieties to toric Varieties}, preprint 1996.


\bibitem{GL2} N. Gonciulea and V. Laksmibai, 
{\em Schubert Varieties, Toric Varieties, and Ladder 
determinantal Varieties}, preprint 1996. 

\bibitem{GMS} B.R. Greene, D.R. Morrison, A. Strominger,
      {\em Black hole condensation and the unification of string vacua},
   Nucl. Phys. {\bf B451} (1995) 109, hep-th/9504145

\bibitem{DM1} D. Morrison, {\em Through the Looking Glass}, alg-geom/9705028.


\bibitem{K3} B. Kim, {\em Quantum Hyperplane Section 
Theorem}, preprint, 1997.  

\bibitem{li-tian} J. Li and G. Tian, {\em Virtual moduli cycles
 and Gromov-Witten invariants}, 
preprint 1996.

\bibitem{LT} A. Libgober and J. Teitelbaum, {\em Lines on Calabi-Yau 
Complete Intersections, mirror symmetry and Picard-Fuchs equations}, 
Duke Math. J., {\bf 69}(1), 29-39 {\em Int. Math. Res. Notes} (1993).  

\bibitem{LS} M. Lynker, R. Schimmrigk, {\em  
 Conifold Transitions and Mirror Symmetries},  
Nucl.Phys. B484 (1997) 562-582,
(hep-th/9511058 )

\bibitem{R} M. Reid, {\em The moduli space of $3$-folds with $K =0$ 
may nevertheless be irreducible}, Math. Ann. {\bf 278} (1987), 329-334. 



\bibitem{Straten} D. van Straten, {\em  A quintic hypersurface 
in ${\bf P}^4$ with $130$ nodes}, Topology, {\bf 32}, no.4   
(1993),857-864.  

\bibitem{S} A. Strominger, {\em  Massless black hole condensation
    and the unification of string vacua},
   Nucl. Phys. {\bf B451} (1995) 96, hep-th/9504090


\bibitem{St} B. Sturmfels, {\em Gr\"obner Bases and Convex Polytopes}, 
Univ. Lect. Notes, vol. 8, AMS, 1996. 

\bibitem{tjotta} E. Tj{\o}tta, {\em Rational curves on the space of determinantal nets of conics}, { Thesis}, Bergen, (1997).
\end{thebibliography}
\end{document}